\documentclass[11pt]{article}
\usepackage{epsfig}

\def\be{\begin{equation}}
\def\ee{\end{equation}}

\def\nn{\nonumber}
\def\<{\langle}
\def\>{\rangle}
\def\lb{\label}

\def\R{{\bf R}}

\def\N{{\bf N}}

\def\diag{{\rm diag}}

\def\E{{\cal E}}

\def\hb{\vrule height0.18cm width0.14cm $\,$}

\title{Index theory for linear self-adjoint operator equations and
nontrivial solutions for asymptotically linear operator equations}
\author{Yujun Dong\thanks{Partially supported by
the National Natural Science Foundation of China}\\
  Department of
Mathematics, Nanjing Normal University, Nanjing,\\ Jiangsu 210097,
P. R. China\\E-mail: yjdong@eyou.com}
\date{}
\begin{document}

\maketitle

\begin{abstract}
{\it We will first establish an index theory for linear
self-adjoint operator equations. And then with the help of this
index theory we will discuss existence and multiplicity of
solutions for asymptotically linear operator equations by making
use of the dual variational methods and  Morse theory. Finally,
some interesting examples concerning second order Hamiltonian
systems, first order Hamiltonian systems and elliptical partial
differential equations will be presented to illustrate our
results.}
\end{abstract}
\bigskip

.

\renewcommand{\theequation}{\thesection.\arabic{equation}}

\setcounter{section}{0}
\setcounter{equation}{0}

\section{Introduction }

\vskip2mm Let $X$ be an infinite-dimensional separable Hilbert
space with inner product $(\cdot,\cdot)$ and norm $||\cdot||$. Let
$Y\subset X$ be a Banach space with norm $||\cdot||_Y$, and the
embedding $Y\hookrightarrow X$ is compact.    Let $A:Y\to X$ be
continuous, selfadjoint, i.e. $(Ax,y)=(x,Ay)$ for any $x,y\in Y$,
$\Im(A)$ is a closed subspace of $X$ and, $\Im(A)\oplus\ker(A)=X$.
In this paper for any $B\in {\cal L}_s(X)$ we first discuss the
classification theory for
 \be
 Ax+Bx=0\lb{1.1}
 \ee
and then discuss solvability of
 \be
 Ax+\Phi'(x)=0\lb{1.2}
 \ee
where $\Phi:X\to \R$ is differentiable. The main results are as
follows.

{\bf Definition 1.1} For any $B\in {\cal L}_s(X)$, we define
  $$
  \nu_A(B)=\dim \ker(A+B).
  $$

$\nu_A(B)$ is called nullity of $B$. It will be proved in Lemma
2.1 next section that the nullity $\nu_A(B)$ is finite.

{\bf Definition 1.2} For any $B_1,B_2\in {\cal L}_s(X)$ with
$B_1<B_2$, we define
  $$
  I_A(B_1,B_2)=\sum_{\lambda\in [0,1)}\nu_A((1-\lambda)B_1+\lambda B_2);
  $$
and for any $B_1,B_2\in {\cal L}_s(X)$ we define
  $$
  I_A(B_1,B_2)=I_A(B_1,k \hbox{id})-I_A(B_2,k \hbox{id})
  $$
where $\hbox{id}:X\to X$ is the identity map and $k$ id$>B_1,k
\hbox{id}>B_2$ for some real number $k>0$.

\vskip 2mm Let $B_0\in{\cal L}_s(X)$ be fixed and let $i_A(B_0)$
be a prescribed integer associated with $B_0$.

{\bf Definition 1.3} For any  $B\in {\cal L}_s(X)$ we define
   $$
   i_A(B)=i_A(B_0)+I_A(B_0,B).
   $$

\vskip2mm We call $i_A(B)$ index of $B$ and $i_A(B_0)$ is called
initial index. Generally, the index $i_A(B)$ depends also on $B_0$
and the initial index. For some well-known precise operators, we
can give the initial index a special value, so that the index
becomes natural. This will be done in the subsequent sections. The
following proposition is also concerned with a precise example.

{\bf Proposition 1.4} If there exists $B_0\in {\cal L}_s(X)$ such
that $\sum_{\lambda<0}\nu_A(B_0+\lambda id)<+\infty$, we will choose
this integer for $i_A(B_0)$. Then the index defined by Definition
1.3 satisfies
  $$
  i_A(B)=\sum_{\lambda<0}\nu_A(B+\lambda id).
  $$

\vskip2mm For index and nullity defined before we have the following properties.

{\bf Proposition 1.5} (i)For any $B,B_1, B_2\in {\cal L}_s(X)$,
$I_A(B_1,B_2)$ and $i_A(B)$ are well-defined and finite;

  (ii) For any $B_1,B_2, B_3\in {\cal L}_s(X)$,
   $I_A(B_1,B_2)+I_A(B_2,B_3)=I_A(B_1,B_3)$;

  (iii) for any $B_1,B_2\in {\cal L}_s(X)$, $I_A(B_1,B_2)=i_A(B_2)-i_A(B_1)$;

  (iv) for any $B_1,B_2\in {\cal L}_s(X)$ with $B_1<B_2$, $ \nu_A(B_1)+i_A(B_1)\leq i_A(B_2)$.

\vskip3mm The following theorems are concerned with solvability of
equation (1.2).

{\bf Theorem 1.6} Assume that $B:X\to {\cal L}_s(X)$ satisfies

 (1) $\Phi'(x)-B(x)x$ is bounded;

 (2) there exist $B_1,B_2\in {\cal L}_s(X)$ satisfying $i_A(B_1)=i_A(B_2), \nu_A(B_2)=0$
 and for any $x\in X$
 $$
 B_1\leq B(x)\leq B_2.
 $$
Then (1.2) has at least one solution.

{\bf Theorem 1.7} Assume that $\Phi''(x)$ is continuous and
bounded,
  $\Phi'(\theta)=\theta$, and

(1) there exist  $B_1,B_2\in {\cal L}_s(X)$ satisfying
 $i_A(B_1)=i_A(B_2), \nu_A(B_2)=0$ such that
   $$
   B_1\leq \Phi''(x)\leq B_2
   $$
for any $x\in X$ with $||x||\geq r>0$;

   (2) with $B_0:=\Phi''(\theta)$, we have
   $$
   i_A(B_1)\notin[i_A(B_0), i_A(B_0)+\nu_A(B_0)].
   $$

Then (1.2) has at least one nontrivial solution $x_0$.
Moreover, under the further assumption that

   (3)$0=\nu_A(B_0)$ and $|i_A(B_1)-i_A(B_0)|\geq \nu_A(\Phi''(x_0))$,

equation (1.2)  has at least two nontrivial solutions.

{\bf Theorem 1.8} Assume that

(1) $\Phi\in C^2(X,\R)$ and there exist $B_1, B_2\in {\cal
L}_s(X)$ with $\nu_A(B_1)=0$ such that
   $$
   B_1\leq \Phi''(x)\leq B_2   \forall x\in X;
   $$

 (2)there exists $B_3\in{\cal L}_s(X)$ with $B_1<B_3$
 and $i_A(B_1)=i_A(B_3), \nu_A(B_3)=0$ such that
  $$
  \Phi(x)\leq {1\over 2}(B_3x,x)+c \forall x\in X;
   $$
  (3) $\Phi'(\theta)=\theta,
  \Phi''(\theta)>B_1,\nu_A(\Phi''(\theta))=0$
  and
  $i_A(\Phi''(\theta))>i_A(B_1)$.

  Then (1.2) has two distinct nontrivial solutions.

{\bf Theorem 1.9} Assume that

  (1) there exist  $B_1,B_2\in {\cal L}_s(X)$ satisfying
 $B_1\leq B_2$ and
 $i_A(B_1)+\nu_A(B_1)=i_A(B_2), \nu_A(B_2)=0$ such that $\Phi(x)-{1\over
 2}(B_1x,x)$ is convex and
   $$
    \Phi(x)\leq {1\over 2}(B_2x,x)+c \forall x\in X.
   $$
Then (1.2) has a solution.

   Moreover, if we further assume that

   (2) $\Phi'(\theta)=\theta, \Phi(\theta)=0$ and there exists $B_0\in {\cal L}_s(X)$ satisfying
$B_0\geq B_1$ and
   $$
   i_A(B_0)>i_A(B_1)+\nu_A(B_1).
   $$
Then (1.2) has at least one nontrivial solution.

\vskip4mm In 1980, H. Amann and E. Zehnder[1] discussed equation
(1.2) under the assumption that $A:$dom$(A)\subseteq X\to X$ is a
unbounded selfadjoint operator. By the saddle point reduction
methods they obtained some existence results for nontrivial
solutions. They also discussed semilinear elliptic boundary value
problems, periodic solutions of a semilinear wave equations, and
periodic solutions of Hamiltonian systems as special cases of the
abstract equation. In 1981, K. C. chang[2] extended their results
by a simpler and unified approach. Especially, Chang obtained an
existence result for three distinct solutions. Theorem 1.8 is
motivated by his this result. And theorem 1.7 comes from his
another result in his excellent book[3]. Chang [4]also discussed
equation (1.2) by assuming that $A\in {\cal L}_s(X)$ and $\Phi'$
is compact. This framework can be used to discuss elliptic partial
differential equations. In 1990, I Ekeland[5] discussed
solvability of equation (1.2) by the dual variational methods and
convex analysis theory. He assumed that $A:X\to X^{\ast}$ is
closed and seladjoint. As applications he mainly focussed on
second order and first order Hamiltonian systems satisfying
various boundary conditions. Our theorem 1.9 generalizes his
results. Some other special equations were also mentioned in the
end of Chapter III.

We would like to stress that our equation (1.2) with the
assumptions on the operator $A$ supplies a new
framework for some special equations. It only requires that the
operator $A$ has finite multiplicities for every eigenvalue.
Most operators listed by Ekeland[5] have this property. Although
this framework can not be used for studying wave equations.
However, it can be used to study seconder order Hamiltonian
systems, first order Hamiltonian systems as well as elliptic
equations. As one can find in sections 3,4 and 5 we will obtain
some new results. One can also find that  the assumptions on the
operator $A$ also make us possible to establish an index theory
for equation (1.1) by the dual variational methods first.
Definitions 1.1,1.2 and 1.3, and propositions 1.4 and 1.5 are
concerned with this index theory. Theorems 1.6,1,7,1,8 and 1.9 can
be regarded as applications of this index theory. Just because of
the usage of this index theory in the assumptions of these theorems
we can get some new results. Note that some special cases 
have been discussed by the author in [6-8].

As far as the author knows, an index theory for convex linear
Hamiltonian systems was established first by I. Ekeland[9] in
1984. By the works[10-13] of Conley, Zehnder and Long, an index
theory for symplectic paths was introduced. These index theories
have important applications[14-22]. One can refers to the two
excellent books[5,23] for systematical treatments. In [24,25] Long
and Zhu defined spectral flows for paths of linear operators and relative
Morse index, and redefined Maslov index for symplectic paths. Our
concept of relative Morse index comes from their papers with some
modifications. Note that in definition 1.2 the  relative Morse
index depends only on nullities. There are also other
contributions on index theory. For example, in 1994, S. E.
Cappell, R. lee and E. Y. Miller[26] introduced three equivalent
definitions for Maslov index. For related topics one can refer to
references[27-29].

The paper is organized as follows. In section 2, we will introduce
an Ekeland type of index theory and prove these results. Sections
3-5 will devote to applications in some special cases of equation
(1.2). Precisely, in section 3 we will first discuss second order
Hamiltonian systems. Then in section 4 we will discuss first order
Hamiltonian systems. Finally , in section 5 we will discuss
elliptic partial differential equations.

\setcounter{equation}{0}
\section{Ekeland type of index theory and proofs of main results}

In his excellent book[5] Ekeland introduced an index theory for
convex linear Hamiltonian systems by dual variational methods. He
also mentioned that by Lasry's tricks some non-convex Hamiltonian
systems could be changed into convex systems and hence could be
discussed also by dual variational methods. In this section we
will make use of his ideas to establish an index theory for linear
system (1.1) first. And then we will prove the main results in the
previous section.

Let $X,Y$ and $A:Y\to X$ be defined as before.

{\bf Lemma 2.1}. For any $B\in {\cal L}_s(X)$, we have that
$A+B:Y\to X$ is continuous, $\ker(A+B)$ is finitely dimensional,
$\Im(A+B)$ is closed and
  $$
  X=\ker (A+B)\oplus \Im(A+B).
  $$

{\bf Proof} Because $\Im(A)\oplus\ker(A)=X$, $Y\bigcap\Im(A)$ is
also a Banach space with the norm $||\cdot||_Y$. So
$A_0:=A|_{\Im(A)\bigcap Y}:\Im(A)\bigcap Y\to \Im(A)$ is
invertible and the inverse $A_0^{-1}:\Im(A)\to \Im(A)\bigcap
Y\hookrightarrow\Im(A)$ is compact and self-adjoint. By the
spectral theory there is a basis $\{e_j\}$ of $\Im(A)$ and a
nonzero sequence $\lambda_j\to 0$ in $\R$ such that:
  \begin{eqnarray}
  &&(e_i,e_j)=\delta_{ij}\lb{2.1}\\
  &&(A_0^{-1}e_j,u)=(\lambda_je_j,u),\,\,\forall\,\,
  u\in \Im(A).\lb{2.2}
  \end{eqnarray}
For any $j\in\N$ we also have $\dim
\ker(A_0^{-1}-\lambda_j)<+\infty$ and $Ae_j={1\over
\lambda_j}e_j$. Fix $k\in\R\setminus\{0\}$ such that ${1\over
\lambda_j}+k\neq 0$ for any $j$. Then if
$u=\sum_{j=1}^{\infty}c_je_j+e_0\in X$ with $e_0\in \ker(A)$, we
have $(A+k id)(\sum_{j=1}^{\infty}({1\over
\lambda_j}+k)^{-1}c_je_j+{1\over
k}e_0)=\sum_{j=1}^{\infty}c_je_j+e_0$, and $\Im(A+k id)=X$. And
for any $B\in {\cal L}_s(X)$, there exists a constant $k$ such
that $ker (A+k id)=\{\theta\}$ and $B-k id>0$. Under a new inner
product defined by $(x,y)_1:=((B-k id)^{-1}x,y)$, $X$ is also a
Hilbert space and $(B-k id)(A+k id)^{-1}:X\to X$ is compact and
selfadjoint. By the spectral theory again, there exists a basis
$\{\zeta_j\}$ of $X$  and a nonzero sequence $\mu_j\to 0$ such
that $(B- k id)(A+k id)^{-1}\zeta_j=\mu_j \zeta_j$. Then $(A+k
id)^{-1}\zeta_j:=\xi_j$  is a basis of $Y$. This means that
$(A+B)\xi_j=(A+k id+B-k id)\xi_j=(1+\mu_j)(A+k id)\xi_j$. So
$\Im(A+B)=\{\sum_{\lambda_j\neq-1}c_j\zeta_j|\sum c_j^2<+\infty\}$
and $\ker(A+B)=\{\sum_{\mu_j=-1} c_j(B-k id)^{-1}\zeta_j\}$ is
finite dimensional. Because $(B-k id)^{-1}>0$ and $((B-k
id)^{-1}\zeta_j,\zeta_i)=0$ for $\mu_j\neq\mu_i$, the projection
of span$\{(B-k id)^{-1}\zeta_j\}_{\mu_j\neq-1}$ to
span$\{\zeta_j\}_{\mu_j\neq-1}$ is
span$\{\zeta_j\}_{\mu_j\neq-1}$. This completes the proof.
\hfill\hb

\vskip4mm From this lemma  for given $B_0\in {\cal L}_s(X)$ we
know that $\Lambda:=(A+B_0)|_{Im(A+B_0)}:Im(A+B_0)\bigcap Y\to
Im(A+B_0)$ is invertible and the inverse
$\Lambda^{-1}:Im(A+B_0)\to X$ is compact. For any $B\in {\cal
L}_s(X)$ with $B-B_0\geq \epsilon id$ for some constant $\epsilon>0$
we define a bilinear form:
  \be
  \psi_{A,B|B_0}(x,y)=(\Lambda^{-1}x,y)+((B-B_0)^{-1}x,y).\lb{2.3}
  \ee
Note that under the inner product $((B-B_0)^{-1}x,y)$, $Im(A+B_0)$
is a Hilbert space, and $(B-B_0)\Lambda^{-1}$ is self-adjoint and
compact. So there exists a basis $\{x_j\}$ of $Im(A+B_0)$
satisfying $(\Lambda^{-1}x_j,x)=\lambda_j((B-B_0)^{-1}x_j,x)$ 
for every $x\in Im(A+B_0)$,
$((B-B_0)^{-1}x_i,x_j)=\delta_{ij}$ and $\lambda_j\to
0$ . Therefore, for any $x=\sum
c_jx_j$ satisfying $\sum_{j=1}^{\infty}c_j^2<+\infty$, we have
  \be
  \psi_{A,B|B_0}(x,x)=\sum_{j=1}^{\infty}(1+\lambda_j)c_j^2.\lb{2.4}
  \ee
Define
 \begin{eqnarray}
 && E_A^+(B|B_0):=\{\sum_{j=1}^{\infty}\xi_je_j|\xi_j=0\,\,\hbox{if}\,\,1+\lambda_j\leq 0\},\nn\\
 && E_A^0(B|B_0):=\{\sum_{j=1}^{\infty}\xi_je_j|\xi_j=0\,\,\hbox{if}\,\,1+\lambda_j\neq 0\},\nn\\
 && E_A^+(B|B_0):=\{\sum_{j=1}^{\infty}\xi_je_j|\xi_j=0\,\,\hbox{if}\,\,1+\lambda_j\geq
 0\}.\nn
 \end{eqnarray}

Then $E_A^0(B|B_0)$ and $E_A^-(B|B_0)$ are finitely dimensional.

{\bf Definition 2.3} For any $B\in {\cal L}_s(X)$ with
$B-B_0>\epsilon id$ we define
 $$
 i_A(B|B_0)=\dim E_A^-(B|B_0), \nu_A(B|B_0)=\dim E_A^0(B|B_0).
 $$

\vskip2mm

This index $(i_A(B|B_0), \nu_A(B|B_0))$ is also a kind of relative
index. Different from $I_A(B_0,B)$ with fixed $B_0$, $i_A(B|B_0)$
is only defined for $B>B_0$. The following theorem lists some
properties concerning these index theories.

 {\bf Theorem 2.4} (1) For any $B>B_0$ we have
   $$
   \nu_A(B|B_0)=\nu_A(B)
   $$
  (2) Assume $B_2>B_1>B_0$, then
   $$
   i_A(B_2|B_0)\geq \nu_A(B_1|B_0)+i_A(B_1|B_0)
   $$
   (3) For any $B_2>B_1>B_0$, we have
  $$
  i_A(B_2|B_0)-i_A(B_1|B_0)=I_A(B_1,B_2).
   $$
   (4) For $B> B_0$, we have
   $i_A(B|B_0)=I_A(B_0,B)-\nu_A(B_0)$.

   (5) With the norm $||u||_{\pm}:=(\pm((\Lambda^{-1}x,x)+((B-B_0)^{-1}x,x)))^{1\over 2}$,
   $E_A^{\pm}(B|B_0)$ are Banach spaces respectively.

   (6) Assume $B_0<B_1<B_2$ with $i_A(B_1)=i_A(B_2)$ and $\nu_A(B_2)=0$. Then
   $$
   X=E^-_A(B_1|B_0)\bigoplus
   E^+_A(B_2|B_0)
   $$

{\bf Proof}. (1) By definition, for any $u\in E_A^0(B|B_0), v\in
Im(\Lambda)$ we have
 $$
 q_{A,B|B_0}(u,v)=0.
 $$
From lemma 2.1, there exists $\xi_u\in\ker(\Lambda)$ such that
 $$
 \Lambda^{-1}u+(B-B_0)^{-1}u=\xi_u
 $$
and
 $$
 \xi_{c_1u_1+c_2u_2}=c_1\xi_{u_1}+c_2\xi_{u_2}.
 $$
Set $x:=\Lambda^{-1}u-\xi_u$. Then $u=\Lambda(x-\xi_u)=\Lambda x$
and
  $$
  Ax+Bx=0.
  $$
So $E_A^0(B|B_0)\cong\ker(A+B)$.

 (2)We only sketch the proof here. We first prove that
 $i_A(B|B_0)$ is a kind of Morse index: for any subspace $X_1$ of
 $\Im(A+B_0)$ satisfying $\psi_{A,B|B_0}(x,x)<0$
 on $X_1\backslash\{\theta\}$, we have that $\dim(X_1)\leq
 i_A(B|B_0)$. Then we check : for any $x\in E_A^-(B_1|B_0)\bigoplus
 E_A^0(B_1|B_0)\backslash\{\theta\}$, $\psi_{A,B_2|B_0}(x,x)<0$.
 Note that in the assumption we only need suppose that
 $((B_2-B_1)x,x)>0$ for any $x\in E_A^+(B_1|B_0)\setminus\{\theta\}$
 in stead of $B_2>B_1$. So in (iv) of proposition 1.5 we can
 assume that $((B_2-B_1)x,x)>0$ for any $x\in \ker(A+B_1)\setminus\{\theta\}$
 in stead of $B_2>B_1$.

 (3) Write $i(s)=i_A((B_1+s(B_2-B_1))|B_0),
 \nu(s)=\nu_A((B_1+s(B_2-B_1))|B_0)$ for $s\in [0,1]$. From Lemma
 2.6, $0\leq i(0)\leq i(1)<+\infty$ and there are only finite many
 $s\in[0,1]$ such that $\nu(s)\neq 0$. For $s\in [0,1]$ with
 $\nu(s)=0$, $i(s)$ is continuous. And for $s\in[0,1]$ with
 $\nu(s)\neq 0$ we have $i(s+0)=i(s-0)=\nu(s)$.

(5) From (2.4), for any $x\in E_A^+(B|B_0)$, we have
$x=\sum_{1+\lambda_j>0}c_jx_j$ with
$\sum_{1+\lambda_j>0}c_j^2<\infty$. So
$||x||_+=(\sum_{1+\lambda_j>0}(1+\lambda_j)c_j^2)^{1\over 2}$ is a
norm.

(6)If $u\in E^-_A(B_1|B_0)\backslash\{\theta\}$, then $\psi_{A,
B_1|B_0}(u,u)<0$ and $\psi_{A, B_2|B_0}(u,u)\leq\psi_{A,
B_1|B_0}(u,u)<0$, and $u\notin E^+_A(B_2|B_0)$. So
$E^-_A(B_1|B_0)\bigcap E^+_A(B_2|B_0)=\{\theta\}$ and we need only
prove $X=E^-_A(B_1|B_0)+ E^+_A(B_2|B_0)$. In fact, by definition
we have $X=E^-_A(B_2|B_0)\bigoplus E^+_A(B_2|B_0)$, and
$i_A(B_2|B_0)=\dim E^-_A(B_2|B_0)<\infty$. Let
$\{e_j\}_{j=1}^{\gamma}$ be a basis of $E^-_A(B_1|B_0)$ where
$\gamma:=i_A(B_1|B_0)$. We have a decomposition $e_j=e_j^-+e_j^+$
with $e_j^-\in E^-_A(B_2|B_0)$ and $e_j^+\in E^+_A(B_2|B_0)$. If
$\sum_{j=1}^{\gamma}\alpha_je_j^-=0$, then ${\bar
x}:=\sum_{j=1}^{\gamma}\alpha_je_j=\sum_{j=1}^{\gamma}\alpha_je_j^+\in
E^+_A(B_2|B_0)$, and ${\bar x}\in E^-_A(B_1|B_0)$. So ${\bar
x}=\theta$ and $\alpha_j=0,j=1,2,\cdots, \gamma$. Hence
$\{e_j^-\}_{j=1}^{\gamma}$ is linear independent. Since $\dim
E^-_A(B_2|B_0)=i_A(B_2|B_0)=i_A(B_1|B_0)=\gamma$, 
$\{e^-_j\}_{j=1}^{\gamma}$ is a basis of $E^-_A(B_2|B_0)$. If
$u\in X, u=u^-+u^+$ with $u^-\in E^-_A(B_2|B_0)$ and $u^+\in
E^+_A(B_2|B_0)$, then $u^-=\sum_{j=1}^{\gamma}\beta_je^-_j$.
So $u=\sum_{j=1}^{\gamma}\beta_je_j+(u^+-\sum_{j=1}^{\gamma}\beta_je_j^+):=u_1+u_2$,
and $u_1\in E^-_A(B_1|B_0)$ and $u_2\in E^+_A(B_2|B_0)$. \hfill\hb

\vskip2mm
 {\bf Proof of Proposition 1.5}. We only prove(i). For
any $B_1<B_2$, by (iii) of theorem 2.4, $I_A(B_1,B_2)$ is finite
and $I_A(B_1,B_2)+I_A(B_2,B_3)=I_A(B_1,B_3)$ if we further assume
$B_2<B_3$. So if $B_1,B_2\in {\cal L}_s(X)$ without any
restriction, there exists $\lambda_0<0$ such that
$\nu_A(B_0+\lambda)=0$ for any $\lambda\leq\lambda_0$. It follows
that $I_A(B_1,k \hbox{id})-I_A(B_2,k \hbox{id})=I_A(B_1,k_1
\hbox{id})-I_A(B_2,k_1 \hbox{id})$ for any $k,k_1\in\R$. So the
relative Morse index and hence the index $i_A(B)$ are finite and
well-defined. \hfill\hb

\vskip3mm

{\bf Proof of Proposition 1.4}. From the additive property,
  $$
  I_A(B_0+\lambda,B_0)=\sum_{\lambda_0\leq\lambda<0}\nu_A(B_0+\lambda)=i_A(B_0).
  $$
So $i_A(B_0+\lambda)=0$ if $\lambda\leq\lambda_0$. For any $B\in
{\cal L}_s(X)$ there exists $\lambda_1<0$ with
$\lambda_1+B<B_0+\lambda_0$. By the monotonicity of indices we
have $i_A(B+\lambda)=0$ and $\nu_A(B+\lambda)=0$ for
$\lambda\leq\lambda_1$.So $I_A(B+\lambda_1,B_0+\lambda_0)\leq
I_A(B_0+\lambda,B_0+\lambda_0)$=0 where $\lambda<0$ is large
enough. So
$i_A(B+\lambda_1)=i_A(B_0+\lambda_0)-I_A(B+\lambda_1,B_0+\lambda_0)=0$.
There exists $\lambda_2<0$ such that $\nu_A(B+\lambda_2)=0$. And
hence,
  $$
  i_A(B)=I_A(B+\lambda_1,B)+i_A(B+\lambda_1)
    =\sum_{\lambda_1\leq\lambda<0}\nu_A(B+\lambda)
    =\sum_{\lambda<0}\nu_A(B+\lambda).
  $$            \hfill\hb

\vskip2mm

{\bf Proof of Theorem 1.6} Choose $k\in\R$ with $\nu_A(k
\hbox{id})=0$. From lemma 2.1, $(A+k \hbox{id})^{-1}$ is compact.
And we need only verify that solutions of the following equations
are a priori bounded:
 $$
 x-k(A+k \hbox{id})^{-1}x+(A+k \hbox{id})^{-1}(\lambda B_1x+(1-\lambda)\Phi'(x))=0
 $$
If not, there exist $\{x_n\}\subset X$ with $||x_n||\to +\infty$
and $\lambda_n\in[0,1]$ such that
 $$
 x_n-k(A+k \hbox{id})^{-1}x_n+(A+k \hbox{id})^{-1}(\lambda
 B_1x_n+(1-\lambda)\Phi'(x_n))=0.
 $$
Set $y_n=x_n/||x_n||$ and $h(x)=\Phi'(x)-B(x)x$. Then
 $$
 y_n-k(A+k \hbox{id})^{-1}y_n+(A+k \hbox{id})^{-1}(\lambda
 B_1y_n+(1-\lambda)(B(x_n)y_n+  ||x_n||^{-1}h(x_n)))=0.
 $$
From the bounded-ness of $B(x)$, for any $y\in X$,
$B(x_n)y\rightharpoonup y_1$. We define ${\bar B}y=y_1$. Then
${\bar B}\in{\cal L}_s(X)$ and $B_1\leq {\bar B}\leq B_2$. By the
compactness of $(A+k \hbox{id})^{-1}$ and the above equation, we
have $y_n\to y_0$ and $B(x_n)(y_n-y_0)\to 0$. We also assume that
$\lambda_n\to \lambda_0$. Taking the limit we have
  $$
  Ay_0+(\lambda_0B_1+(1-\lambda_0){\bar B})y_0=0.
  $$
But $B_1\leq B_3:=\lambda_0B_1+(1-\lambda_0){\bar B}\leq B_2$
leads to $\nu_A(B_3)=0$. This is a contradiction to the fact that
$y_0$ is a nontrivial solution.\hfill\hb

\vskip4mm
In the following we will prove Theorem 1.7. To do this
we need a lemma, which comes from [6, Chapter II. Theorem 5.1, 5.2 and Corollary 5.2]. Note that for
any $B\in {\cal L}_s(X)$, $m^-(B)$ denotes the multiplicity of the
negative eigenvalues of $B$ and $m^0(B)$ denotes the multiplicity
of zero eigenvalues of $B$.

{\bf Lemma 2.5}. Assume $f\in C^2(X,R)$ satisfies the (PS)
condition, $f'(\theta)=\theta$, and there is a positive integer
$\gamma$  such that
$\gamma\notin[m^-(f''(\theta)),m^0(f''(\theta))+m^-(f''(\theta))]$
and $H_q(X,f_a;\R)=\delta_{q\gamma}\R$ for some regular value
$a<0$. Then $f$ has a critical point $p_0\neq \theta$ with
$C_{\gamma}(f,p_0)\neq0$. Moreover, if $\theta$ is a
non-degenerate critical point, and
$m^0(f''(p_0))\leq|\gamma-m^-(f''(\theta))|$, then $f$ has another
critical point $p_1\neq p_0,\theta$.

\vskip3mm We now begin to prove Theorem 1.7.  From assumption (1)
and that $\Phi''(x)$ is bounded we can choose $k_1,k\in\R$ such
that $\nu_A(k id)=0=\nu_A(k_1 id), B_1(t)-k id\geq id$ and
   \begin{eqnarray}
  &&k_1 id\geq N''(x)\geq id\,\,\forall x\lb{2.5}\\
  &&B_2-k id\geq N''(x)\geq B_1-k id\,\,\hbox{for}\,\,|x|\geq
  r,\lb{2.6}
  \end{eqnarray}
where  $N(x)=\Phi(x)-{1\over 2}k(x,x)$. By the (iii) of
Proposition 1.5, we may assume $\nu_A(k id)=0=\nu_A(k_1 id)$. Let
$\Lambda u:=Au+ku$ and consider the functional
 \be
 \psi(u)={1\over
 2}(\Lambda^{-1}u,u)+N^{\ast}(u)\,\,\forall u\in
 X.\lb{2.7}
 \ee
We have the following proposition.

\vskip2mm
 {\bf Proposition 2.6}  Under the assumption (i) in
theorem 1.7, the functional $\psi$ defined by (2.7) satisfies the
$(PS)$ condition.

{\bf Proof} Assume $\{u_j\}\subset X$ such that $\psi(u_j)$ is
bounded and $\psi'(u_j)\to \theta$ in $X$. If $||u_j||_X$ is
bounded, then there exists a subsequence $u_{j_k}\rightharpoonup
u_0$ in $X$, and $\Lambda^{-1}u_{j_k}\to \Lambda^{-1}u_0$. From
the following (6.9), we have
$\Lambda^{-1}u_{j_k}+{N^{\ast}}'(u_{j_k})=\psi'(u_{j_k})$, and
${N^{\ast}}'(u_{j_k})=\psi'(u_{j_k})-\Lambda^{-1}u_{j_k}\to
-\Lambda^{-1}u_0$ in $X$. By the Fenchel conjugate formula and [1,
Theorem II.4], $u_{j_k}=N'(\psi'(u_{j_k})-\Lambda^{-1}u_{j_k})\to
N'(-\Lambda^{-1}u_0)$ in $X$ and $\psi$ satisfies the $(PS)$
condition. So in the following we only need to show $\{u_j\}$ is
bounded in $X$.

From $N'(\theta)=\theta$, we have ${N^{\ast}}'(\theta)=\theta$ and
 \be
 (\psi'(u),v)=(\Lambda^{-1}u,v)+({N^{\ast}}'(u),v)\,\,\forall
 v,u\in X.\lb{2.8}
 \ee
Noticing that $\int_0^1{N^{\ast}}''(\theta u_j)d\theta u_j
={N^{\ast}}'(u_j)$, we have
 \be
 \Lambda^{-1}u_j+\int_0^1{N^{\ast}}''(s u_j)ds u_j
 =\psi'(u_j)\to \theta,\,\,\hbox{in}\,\,X.\lb{2.9}
 \ee
If $||u_j||_X$ is not bounded, without loss of generality we
assume $||u_j||_X\to\infty$. Set $x_j=u_j/||u_j||_X$. We also
assume $x_j\rightharpoonup x_0$ in $X$ by going to subsequence if
necessary. And hence $\Lambda^{-1}x_j\to \Lambda^{-1}x_0$ in $X$.
From [5,Propositions II.2.10, I.1.15] and (2.5) we have
${N^{\ast}}''(u^{\ast})=(N''(u))^{-1}$ as $u^{\ast}=N'(u)$, and
  \begin{eqnarray}
  &&id\leq {N^{\ast}}''(x)\leq k_1^{-1}id\,\,\forall\,\,x\in X\lb{2.10}\\
  &&(B_2-k id)^{-1}\leq {N^{\ast}}''(x)\leq (B_1- k id)^{-1}\,\,\forall\,\,x\in
  X\,\,\hbox{with}\,\,||x||\geq r_1.\lb{2.11}
  \end{eqnarray}
For any $\delta\in (0,1)$ fixed, set
 \begin{eqnarray}
C_j&&=\int_0^1{N^{\ast}}''(su_j)ds,
         ||u_j||\geq {r_1\over \delta}\nn\\
       &&=(B_1-k id)^{-1},\,\,\hbox{otherwise}\nn
 \end{eqnarray}
and
 $$
\xi_j=\int_0^1{N^{\ast}}''(s u_j)ds u_j-C_ju_j.
 $$
From assumption(1) and (2.5)(2.9), there exists a constant $c_1>0$
such that
 \begin{eqnarray}
    &&||\xi_j||\leq c_1,\lb{2.12}\\
   && (1-\delta)(B_2-k id)^{-1}+\delta id\leq
  C_j\leq(1-\delta)(B_1-k id)^{-1}+k_1^{-1}\delta id\lb{2.13}\\
  && \Lambda^{-1}u_j+C_ju_j+\xi_j=\psi'(u_j).\lb{2.14}
 \end{eqnarray}
Now by going to subsequences if necessary we may further assume
$C_ju\rightharpoonup C_0 u$ in $X$ for every $u\in X$. And from
(2.9)(2.13)(2.14), for every $\epsilon>0$ we have
  \begin{eqnarray}
  &&(B_2-(k+\epsilon)id)^{-1}\leq
    C_0\leq (B_1-(k-\epsilon)id)^{-1},\nn\\
  &&\Lambda^{-1}x_0+C_0x_0=0\nn
  \end{eqnarray}
Let $\Lambda^{-1}x_0=y_0$ and $B_0=C_0^{-1}+k id$. We have
 \be
 Ay_0+B_0y_0=0.\lb{2.15}
 \ee
We need only show that this is a contradiction. From assumption(1)
and the finiteness of the relative Morse index, for $\epsilon>0$
is small enough, we have $\nu_A(B_1-\epsilon
id)=\nu_A(B_2+\epsilon id)=0$ and $i_A(B_1-\epsilon
id)=i_A(B_2+\epsilon id)$. So that $B_1-\epsilon id\leq B_0\leq
B_2+\epsilon id$ and $\nu_A(B_0)=0$. This is impossible since
$||y_0||_X=1$ and $y_0$ is a nontrivial solution of (2.15). This
contradiction means $||u_j||_X$ is bounded. \hfill\hb

\vskip2mm
{\bf Proof of Theorem 1.7}. From Lemma 2.5 and
Proposition 2.6 it suffices to show that
 \be
 H_q(X,\psi_{-a};\R)\cong\delta_{q\gamma}\R,
 q=0,1,2,\cdots,\lb{2.16}
 \ee
for $a>0$ is large enough, where $\gamma:=i_A(B_1|k id)$. In fact
${N^{\ast}}''(\theta)=(N''(\theta))^{-1}$ and
  \begin{eqnarray}
  (\psi''(\theta)u,u)&&=(\Lambda^{-1}u,u)+({N^{\ast}}''(\theta)u,u)\nn\\
    &&=(\Lambda^{-1}u),u)+((B_0-k id)^{-1}u,u), \,\,\forall u\in X.\nn
  \end{eqnarray}
By definition, $m^-(\psi''(\theta))=i_A(B_0|k id),
m^0(\psi''(\theta))=\nu_A(B_0|k id)$.  $i_A(B_1)\notin
[i_A(B_0),i_A(B_0)+\nu_A(B_0)]$ if and only if $i_A(B_1|k id)\in
[i_A(B_0|k id), i_A(B_0)-i|k id B)+\nu_A(B_0|k id)]$; and
$\nu_A(B_0)=\nu_A(B_0|k id), |i_A(B_1)-i_A(B_0)|=|i_A(B_1|k
id)-i_A(B_0|k id)|$.

We will prove (2.16) in the following two steps.

{\bf Step 1}. For $\epsilon>0$ is sufficiently small, set ${\cal
M}_R:=(E^+_A(B_2+\epsilon id|k id)\bigcap B_R)\bigoplus
E^-_A(B_1-\epsilon id|k id)$, then for $R,a>0$ are large enough we
have
  \be
  H_q(X,\psi_{-a};\R)=H_q({\cal M}_R,{\cal M}_R\cap \psi_{-a}; \R),
  q=0,1,2,\cdots\lb{2.17}.
  \ee
In fact, for any $\epsilon>0$ is small enough we have
$i_A(B_1-\epsilon id)=i_A(B_2+\epsilon id)$ and
$\nu_A(B_2+\epsilon id)=0$. It is easy to see $E^-_A(B_1-\epsilon
id|k id)$ and $E^+_A(B_2+\epsilon id|k id)$ are Banach spaces
under the following norms
  $$
 ||u||_1:=((\Lambda^{-1}u,u)
    +((B_1-k id-\epsilon id)^{-1}u,u))^{1\over 2}.
  $$
and
 $$
   ||u||_2:=((\Lambda^{-1}u,u)
    +((B_2(t)-k id+\epsilon id)^{-1}u,u))^{1\over 2}
 $$
respectively. So for every $u=u_1+u_2\in X$ with $u_1\in
E^-_A(B_1-\epsilon id|k id)$ and $u_2\in E^+_A(B_2+\epsilon id|k
idB)$, from (2.11) we have
 \begin{eqnarray}
  (\psi'(u),u_2-u_1)&&=(\Lambda^{-1}u,u_2-u_1)+({N^{\ast}}'(-u),u_1-u_2)\nn\\
                   &&=-(\Lambda^{-1}u_1,u_1)+(\int_0^1{N^{\ast}}''(-\theta u)d\theta u_1,u_1)\nn\\
                    &&+(\Lambda^{-1}u_2,u_2)-(\int_0^1{N^{\ast}}''(\theta u)d\theta u_2,u_2)\nn\\
                    &&\geq-(\Lambda^{-1}u_1,u_1)+((B_1-k id-\varepsilon id)^{-1}u_1,u_1)\nn\\
                    &&+(\Lambda^{-1}u_2,u_2)+((\varepsilon id
                    +B_2(t)-k id)^{-1}u_2,u_2)-c_2\nn\\
                    &&\geq
                    c_3||u_2||_X^2+c_4||u_1||_X^2-c_2,\nn
 \end{eqnarray}
where $c_2,c_3, c_4>0$ are constants. When $R$ is large enough we
have
   $$
   (\psi'(u),u_2-u_1)>1
   $$
for every $u=u_1+u_2$ with $u_1\in E^-_P(B_1-\epsilon I_{2n}|B),
u_2\in E^+_P(B_2+\epsilon I_{2n}|B)$ and $||u_2||_{L^2}\geq R$, or
$||u_1||_{L^2}\geq R$. For any $u=u_2+u_1\notin{\cal M}_R$, let
$\sigma(t,u)=e^{-t}u_2+e^tu_1, T_u=\ln ||u_2||-\ln R$, and
  \begin{eqnarray}
  \eta(t,u_2+u_1)&&=u_2+u_1, ||u_2||\leq R,\nn\\
                 &&=\sigma(T_ut,u), ||u_2||>R.\nn
  \end{eqnarray}
Then $\eta:[0,1]\times L^2\to L^2$ is continuous and $({\cal
M}_R,{\cal M}_R\bigcap \psi_{-a})$ is a deformation retract of
$(L^2,\psi_{-a})$. Therefore, (2.17) is satisfied.

{\bf Step 2}. For $R,a>0$ are large enough, we have
 $$
 H_q({\cal M}_R,
   {\cal M}_R\bigcap \psi_{-a};\R)\cong \delta_{q \gamma}\R,
   q=0,1,\cdots.
 $$

In fact, we have from (2.10)(2.11) that
  \begin{eqnarray}
    &&N^{\ast}(u)\nn\\
          &&=(\int_0^1\theta d\theta\int_0^1{N^{\ast}}''(\theta su)ds u,
             u)+N^{\ast}(\theta)\nn\\
          &&=\int_{\delta}^1\theta d\theta\int_{\delta}^1{N^{\ast}}''(\theta su)ds u,
             u)+o(1)(u,u)\nn\\
          &&\leq{1\over 2}((B_1-k id-\epsilon
          id)^{-1}u,u)\nn
 \end{eqnarray}
when $||u||\geq r_1/{\delta}^2$. Here in the second equality
$o(1)\to 0$ as $\delta \to 0$. Hence, for every $\epsilon>0$ there
exists a constant $c_5$ such that
  \be
N^{\ast}(u)\leq (B_1-k id-\epsilon
                   id)^{-1}u,u)+c_5,\,\,\forall u\in X.\lb{2.18}
  \ee
Therefore, for any $u=u_1+u_2$ with $u_1\in E^-_A(B_1-\epsilon
id|k id)$ and $u_2\in E^+_A(B_2+\epsilon id|k id)\bigcap B_R$,
from (2.7)(2.18) we have
  $$
  \psi(u)\leq -c_4 ||u_1||_X^2+
  c_6||u_1||_X+c_7
  $$
where $c_6, c_7>0$. And hence,
 $$
 \psi(u)\to-\infty\Longleftrightarrow ||u_1||\to+\infty\,\,\hbox{uniformly
 in}\,\,u_2\in E^+_A(B_2+\epsilon id|k id)\cap B_R.
 $$

So, there exist $T>0, a_1>a_2>T, 0<R_1<R_2<R_0$ such that
 \be
 {\cal N}_{R_2}\subset
  \psi_{-a_1}\cap{\cal M}_{R_0}\subset{\cal N}_{R_1}\subset \psi_{a_2}\cap{\cal
  M}_{R_0},\lb{2.19}
 \ee
where ${\cal N}_R:=(E^+_A(B_2+\epsilon id|k id)\cap
B_{R_0})\oplus(E^-_P(B_1-\epsilon id|k id)\setminus B_R)$ and
$B_R$ denotes the closed neighborhood of the origin with radius $R$
in a Banach space. For any $u\in {\cal
M}_{R_0}\cap(\psi_{-a_2}\setminus \psi_{-a_1})$, since
$\sigma(t,u)=e^{-t}u_2+e^tu_1, \psi(\sigma(t,x))$ is continuous
with respect to $t, \psi(\sigma(0,x))=\psi(u)>-a_1$ and
$\psi(\sigma(t,u))\to-\infty$ as $t\to+\infty$, so there exists a
unique  $t=T_1(u)$ such that $\psi(\sigma(t,u))=-a_1$. Because
 \begin{eqnarray}
 &&{d\over dt}\psi(\sigma(t,u))=<\psi'(\sigma(t,u)),\sigma'(t,u)>\nn\\
   &&=<\psi'(e^{-t}u_2+e^tu_1),-e^{-t}u_2+e^tu_1>\leq-1\nn
   \end{eqnarray}
as $t>0$, by the implicit function theorem, $t=T_1(u)$ is
continuous. Define
 \begin{eqnarray}
 \eta_1(t,u)&&=u,u\in \psi_{-a_1}\cap{\cal M}_{R_0}\nn\\
            &&=\sigma(T_1(u)t,u),u\in{\cal
            M}\cap(\psi_{-a_2}\setminus \psi_{-a_1})\nn
 \end{eqnarray}
then $\eta_1:[0,1]\times \psi_{-a_2}\cap{\cal M}_{R_0}\to
\psi_{-a_2}\cap{\cal M}_{R_0}$ and $(\psi_{-a_1}\cap{\cal
M}_{R_0},\psi_{-a_1}\cap{\cal M}_{R_0})$ is a deformation retract
of  $(\psi_{-a_2}\cap{\cal M}_{R_0},\psi_{-a_1}\cap{\cal
M}_{R_0})$. Hence,
 \be
 H_q(\psi_{-a_2}\cap{\cal
M}_{R_0},\psi_{-a_1}\cap{\cal M}_{R_0})\cong
H_q(\psi_{-a_1}\cap{\cal M}_{R_0},\psi_{-a_1}\cap{\cal
M}_{R_0})\cong 0.\lb{2.20}
 \ee

Recall that for any topological spaces $Z\subseteq Y\subseteq X$,
we have exact sequences
  $$
  H_q(Y,Z)\to H_q(X,Z)\to H_q(X,Y)\to H_{q-1}(Y,Z).
  $$
From (2.19), in order to prove
  \be
  H_q({\cal M}_{R_0}, {\cal M}_{R_0}\cap \psi_{-a_2})\cong H_q({\cal M}, {\cal
  N}_{R_2})\lb{2.21}
  \ee
we only need to prove
  $$
  H_q({\cal M}_{R_0}\cap \psi_{-a_2},{\cal N}_{R_2})\cong
  0.
  $$
And from (2.20), it suffices to show
  \be
  H_q({\cal M}_{R_0}\cap \psi_{-a_1},{\cal N}_{R_2})\cong
  0.\lb{2.22}
  \ee
Let $\eta:[0,1]\times {\cal N}_{R_1}\to {\cal N}_{R_1}$ satisfy
   \begin{eqnarray}
   \eta(t,u^++u^-)&&=u^++u^-\,\,\hbox{if}\,\,||u^-||\geq R_2\nn\\
                  &&=u^++{u^-\over
                  ||u^-||}(tR_2+(1-t)||u^-||)\,\,\hbox{if}\,\,R_2>||u^-||\geq
                  R_1.\nn
   \end{eqnarray}
Set $\tau_1=\eta_1(1,\cdot)$ and
$\xi=\tau_1\circ\eta:[0,1]\times(\psi_{-a_1}\cap{\cal M
})\to\psi_{-a_1}\cap{\cal M }$. Then $({\cal N}_{R_2},{\cal
N}_{R_2})$ is a deformation retract of $(\psi_{-a_1}\cap{\cal M
}_{R_0},{\cal N}_{R_2})$. As a result, (2.22) and (2.21) are
valid. From (2.21) we have
  \begin{eqnarray}
   &&H({\cal M}_{R_0},{\cal M}_{R_0}\cap \psi_{-a_2};\R)\nn\\
    && \cong H_q({\cal M}_{R_0},{\cal N}_{R_2};\R)\nn\\
    &&\cong H_q(E^-_P(B_1-\epsilon I_{2n}|B)\cap B_{R_1},\partial(E^-_P(B_1-\epsilon I_{2n}|B)\cap
            B_{R_1});\R)\nn\\
    &&\cong\delta_{q\gamma}\R, q=0,1,2\cdots.\nn
 \end{eqnarray}
\hfill\hb

{\bf Proof of Theorem 1.8}. Define $\Lambda x=Ax+B_1x,
N(x)=\Phi(x)+{1\over 2}(B_1x,x), N^{\ast}(x)=\sup_{y\in
X}\{(x,y)-N(y)\}$ and
  \be
 \psi(u)={1\over
 2}(\Lambda^{-1}u,u)+N^{\ast}(u)\,\,\,\, u\in
 X.\lb{2.23}
 \ee
From assumption (2), we have
  \be
\psi(u)\geq{1\over 2}[(\Lambda^{-1}u,u)+(({\bar
B_2}-B_1)^{-1}u,u)]-c\,\,\,\,\forall u\in X.\lb{2.24}
  \ee
So $\psi$ is bounded from below. If $\psi(u_j)$ is bounded, we
have $||u_j||_X$ is also bounded. We can assume
$u_j\rightharpoonup u_0$ and $\Lambda^{-1}u_j\to \Lambda^{-1}u_0$.
If $\psi'(u_j)\to 0$, we have
$N^{\ast}(u_j)=\psi'(u_j)-\Lambda^{-1}u_j\to-\Lambda^{-1}u_0$ and
$u_j=N'(\psi'(u_j)-\Lambda^{-1}u_j)\to N'(-\Lambda^{-1}u_0)$ in
$X$. So $\psi$ satisfies the (PS) condition. It is easy to check
that $\psi''(\theta):X\to X$ is invertible;
$m^-(\psi''(\theta))=i_A(B_0|B_1)>0$, so that $\theta$ is not a
minimal point. From a theorem in [4],this complete the
proof.\hfill\hb

{\bf Proof of Theorem 1.9}. Consider the functional defined in
(2.23). At this time its domain is not $X$ but $Im(\Lambda)$. We
also have inequality (2.24) with ${\bar B_2}$ and $X$ instead with
$B_2$ and $Im(\Lambda)$respectively. $\psi$ is bounded from below.
Let $u_n\in Im(\Lambda)$ satisfying $\psi(u_n)\to Inf
\psi>-\infty$. Then $\{u_n\}$ is bounded and we assume that
$u_n\rightharpoonup u_0$ in $Im(\Lambda)$. By the compactness of
$\Lambda_0^{-1}$ and the weakly lower semi-continuity of
$N^{\ast}$, we have $Inf \psi(u)\geq\psi(u_0)$. This means that
$u_0$ is a critical point of $\psi$. A simple calculation shows
that $(\Lambda_0^{-1}u_0+N^{\ast}(u_0),u)=0$ for any $u\in
Im(\Lambda)$. So $x=\Lambda_0^{-1}+x_0$ is a solution of (1.2) for
some $x_0\in \ker(\Lambda)$. When $\Phi'(\theta)=\theta$, then
$\theta$ is a solution of (1.2). We will prove that
$u_0\neq\theta$ under assumption (2). In fact, we have
  $$
\psi(u)\leq{1\over
2}[(\Lambda_0^{-1}u,u)+((B_0-B_1)^{-1}u,u)]\,\,\hbox{as}\,\,u\to\theta.
  $$
The Morse index of the right functional at $u=\theta$ is
$\dim(E_A^-(B_0|B_1))$. So for any $u\in
E_A^-(B_0|B_1)\backslash\{\theta\}$ small enough, we have
$\psi(u)<0=\psi(\theta)$. Hence $u_0\neq\theta$ and
$x=\Lambda_0^{-1}+x_0\neq\theta$. This completes the
proof.\hfill\hb

\setcounter{equation}{0}
\section{Second order Hamiltonian systems}

 In this section we will make use of proposition 1.4 to give some classifications for
 second order Hamiltonian systems.

 {\bf 3.1 Sturm-Liouville BVPs}

 In this subsection we will establish
  a classification theory for the following Lagrangian system
  satisfying Sturm-Liouville BVPs
  \begin{eqnarray}
  &&(\Lambda(t)x')'+B(t)x=0\lb{3.1}\\
  &&x(0)\cos\alpha-\Lambda(0)x'(0)\sin\alpha=0\lb{3.2}\\
  &&x(1)\cos\beta-\Lambda(1)x'(1)\sin\beta=0\lb{3.3}
 \end{eqnarray}
where $\Lambda\in C([0,1];GL_s(\R^n)), B\in
L^{\infty}((0,1);GL_s(\R))$ is positive definite for $t\in [0,1]$ , $0\leq\alpha<\pi$ and $0<\beta\leq\pi$.
Let $X:=L^2((0,1);GL_s(\R^n)), Y=\{x\in
C^1([0,1],\R^n)|(\Lambda(t)x'(t))'\in L^2((0,1);\R^n),x(t)$
satisfies $(3.2)(3.3)\}$. For every $x\in Y$, we define
$||x||_Y:=(\int_0^1|x(t)|^2+|x'(t)|+|(\Lambda(t)x'(t))'|^2dt)^{1\over
2}$, $A:X\to Y$ by $(Ax)(t):=(\Lambda(t)x'(t))'$. Then $X$ is a
separable Hilbert space, $Y$ is a Banach space and the embedding
$Y\hookrightarrow X$ is compact. Define $(Bx)(t)=B(t)x(t)$ for any $x\in X$. 
Then equation (3.1)(3.2)(3.3) is equivalent to equation (1.1). 
In view of the following lemma 3.1, the following problem (3.2)(3.3) and
   $$
   (\Lambda(t)x')'-({\bar \lambda}+1)x=0
   $$
has no nontrivial solutions. From general theory of ordinary
differential equations(c.f. [30, pp407-408] for example), for any
$h\in L^2((0,1);\R^n)$, the following problem (3.2)(3.3) and
 $$
  (\Lambda(t)x')'-({\bar \lambda+1})x=h(t)
  $$
has a unique solution. So from lemma 2.1, $A$ is continuous and
closed and, $\ker(A)\bigoplus\Im(A)=X$.

{\bf Lemma 3.1} There exists $\bar\lambda>0$ such that
  $$
  \int_0^1(\Lambda(t)x'(t))'\cdot
  x(t)dt\leq{\bar\lambda}\int_0^1|x(t)|^2dt \forall x\in Y.
  $$

{\bf Proof}.As $\alpha=0$, we have $x(0)=0$; as $\alpha\neq 0$, we
have $\Lambda(0)x'(0)=x(0)\cot\alpha$. By partial integration we
have
 $$
  \int_0^1(\Lambda(t)x'(t))')\cdot
  x(t)dt=-\int_0^1\Lambda(t)x'(t)\cdot x'(t)dt+\Lambda(1)x'(1)\cdot x(1)-\Lambda(0)x'(0)\cdot
  x(0).
  $$
So we need only prove that: for any given $a>0$, there exists
$\lambda_a>0$ such that
  $$
  \int_0^1|x'(t)|^2dt+\lambda_a\int_0^1|x(t)|^2dt\geq
  a(|x(0)|^2+|x(1)|^2).
  $$
The following trick comes from Professor Eric Sere:
 \begin{eqnarray}
 &&{d\over dt}|x(t)|^2={1\over 2}x'(t)\cdot x(t)\geq -\epsilon
 |x'(t)|^2-{1\over \epsilon}|x(t)|^2,\nn\\
 &&|x(t)|^2\geq-\epsilon\int_0^1|x(t)|^2dt-{1\over
 \epsilon}\int_0^1|x(t)|^2dt+|x(0)|^2,\nn\\
 &&(1+{1\over
 \epsilon})\int_0^1|x(t)|^2dt+\epsilon\int_0^1|x'(t)|^2dt\geq
 |x(0)|^2\nn
 \end{eqnarray}
where $\epsilon>0$ is a constant. This completes the proof.
\hfill\hb

Lemma 3.1 shows that $\ker(A-\lambda I_n)=\{\theta\}$ for
$\lambda>{\bar \lambda}$. In view of proposition 1.4 we can
give the following definition.

 {\bf Definition 3.2} For any $B\in L^{\infty}((0,1);GL_s(\R^n))$,
 we define
  \begin{eqnarray}
  &&\nu^s_{\Lambda,\alpha,\beta}(B):=\dim \ker(A+B),\nn\\
  &&i^s_{\Lambda,\alpha,\beta}(B):=\sum_{\lambda<0}\nu^s_{\Lambda,\alpha,\beta}(B+I_n\lambda).\nn
  \end{eqnarray}

For any $B_1,B_2\in L^{\infty}((0,1);$GL$_s(\R^n))$, we define
$B_1\leq B_2$ if and only if $B_1(t)\leq B_2(t)$ for a.e. $t\in
(0,1)$; and define $B_1<B_2$ if and only if $B_1\leq B_2$ and
$B_1(t)<B_2(t)$ on a subset of $(0,1)$ with positive measure.

{\bf Proposition 3.3} We have the following property:

 (1) For any $B\in L^{\infty}((0,1);GL_s(\R^n))$, we have
  $$
 \nu^s_{\Lambda,\alpha,\beta}(B)\in\{0,1,\cdots,n\}.
  $$

 (2) For any $B_1, B_2\in L^{\infty}((0,1);$GL$_s(\R^n))$
satisfying $B_1< B_2$, we have
$i^s_{\Lambda,\alpha,\beta}(B_1)+\nu^s_{\Lambda,\alpha,\beta}(B_1)\leq
i^s_{\Lambda,\alpha,\beta}(B_2)$.

 {\bf Proof} (1)Let $y(t)=\Lambda(t)x'(t), z=(y,x)$, then
(3.1)-(3.3) has an
 equivalent form:
  \begin{eqnarray}
  &&\dot{z}=J\diag\{\Lambda(t)^{-1},B(t)\}z,\lb{3.4}\\
  &&x(0)\cos\alpha-y(0)\sin\alpha=0,\nn\\
  &&x(1)\cos\beta-y(1)\sin\beta=0.\nn
  \end{eqnarray}
 Let $\gamma(t)$ be the fundamental solution of (3.4). Then
   \begin{eqnarray}
   \ker(A+B)&=&\{z(t)=\gamma(t)c|c\in\R^{2n}, z=(y,x) \hbox{satisfies}
   (3.2)(3.3)\}\nn\\
   &\cong&\{c_1,c_2\in\R^n|c_1\cos\alpha-c_2\sin\alpha=0,(I_n\cos\beta,-I_n\sin\beta)\gamma(1)(c_1,c_2)^{\tau}=0\}\nn\\
   &\cong&\{c\in\R^n|(I_n\cos\beta,-I_n\sin\beta)\gamma(1)(0,c)^{\tau}=0\}\subseteq\R^n,\nn
   \end{eqnarray}
   as $\alpha=0$.

   (2) Follows directly from the (iii) of proposition 1.5 and the proof of (2) of theorem 2.4.  \hfill\hb

We now begin to discuss solvability of the following nonlinear
Hamiltonian systems:
 \begin{eqnarray}
  &&(\Lambda(t)x')'+V'(t,x)=0,\lb{3.5}\\
  &&x(0)\cos\alpha-\Lambda(0)x'(0)\sin\alpha=0\nn\\
  &&x(1)\cos\beta-\Lambda(1)x'(1)\sin\beta=0\nn
  \end{eqnarray}
where $V:[0,1]\times\R^n\to \R^n$ is continuous and $V'(t,x)$ denotes the gradient of $V(t,x)$ with respect to $x$.

Define $\Phi(x)=\int_0^1V(t,x(t))dt$ for every $x\in X$. Then
$\Phi'(x)=V'(\cdot,x(\cdot))$ when $V\in C^1([0,1]\times\R^n,\R)$,
and $\Phi''(x)=V''(\cdot,x(\cdot))$ when $V\in
C^2([0,1]\times\R^n,\R)$. Obviously, equation (3.5)(3.2)(3.3) is equivalent to equation (1.2).
From theorem 1.6, theorem 1.7 and its
proof  we have the following results.

 {\bf Theorem 3.4} Assume that $V\in C^1([0,1]\times \R^n, \R)$ and there exist $B_1, B_2\in
L^{\infty}((0,1);$GL$_s(\R^n))$ with $B_1\leq B_2,
i^s_{\Lambda,\alpha,\beta}(B_1)=i^s_{\Lambda,\alpha,\beta}(B_2)$,
$\nu^s_{\Lambda,\alpha,\beta}(B_2)=0$ and $B\in C([0,1]\times\R^n,
GL_s(\R^n))$ such that
  \begin{eqnarray}
  &&V'(t,x)-B(t,x)x\,\,\hbox{is bounded}\,\,\nn\\
  &&B_1(t)\leq B(t,x)\leq B_2(t), (t,x)\in (0,1)\times\R^n,\,\hbox{with}\,\, |x|\geq r>0\nn.
  \end{eqnarray}
  Then (3.1)(3.2)(3.3) has at least one solution.

{\bf Theorem 3.5}  Assume

(1) $V\in C^2([0,1]\times \R^n, \R), B_1(t)\leq V''(t,x)\leq
B_2(t)$ for $|x|\geq r>0$ with
$i^s_{\Lambda,\alpha,\beta}(B_1)=i^s_{\Lambda,\alpha,\beta}(B_2),\nu^s_{\Lambda,\alpha,\beta}(B_2)=0$.

(2) $V'(t,0)\equiv 0, {\bar B}(t):=V''(t,0)$ and
$i^s_{\Lambda,\alpha,\beta}(B_1)\notin
[i^s_{\Lambda,\alpha,\beta}({\bar B}),
i^s_{\Lambda,\alpha,\beta}({\bar
B})+\nu^s_{\Lambda,\alpha,\beta}({\bar B})$.

Then problem (3.1)(3.2)(3.3) has at least one nontrivial solution.
Moreover, if we further assume

(3)$\nu^s_{\Lambda,\alpha,\beta}({\bar
B})=0,|i^s_{\Lambda,\alpha,\beta}(B_1)-i^s_{\Lambda,\alpha,\beta}({\bar
B})|\geq n$.

Then (3.1)(3.2)(3.3) has two nontrivial solutions.

{\bf Remarks 1}. As $\alpha=0,\beta=\pi,\Lambda(t)\equiv I_n$,
linear system (3.1)(3.2)(3.3) reduces to
 \begin{eqnarray}
  &&x''+B(t)x=0\nn\\
  &&x(0)=0=x(1).\nn
 \end{eqnarray}
 An index theory $(i(B),\nu(B))$ was established in [6](2005) by making use a
 direct variational method. Note that this index theory is a special case of definition 3.2, i.e.,
 $(i(B),\nu(B))=(i^s_{I_n,0,\pi}(B),\nu^s_{I_n,0,\pi}(B))$. The index theory $(i(B),\nu(B))$ was used to discuss
 associated second order nonlinear Hamiltonian systems. Note that
 most of the main results in [6] are covered by theorems 3.4, 3.5.
 For related topics one can refers to [31-33].

   {\bf 2}. As $n=1, \Lambda(t)=1$, equation (3.5) is called Duffing equation as
   usual and can be expressed as
     $$
     x''+f(t,x)=0.
     $$
Many papers devoted to solvability of this Duffing equation
satisfying various boundary conditions(see [34-42] and the
references therein). One can find that the main results of some of
these papers are special cases of theorem 3.4 or the following
theorem 3.10.

   {\bf 3.2 Generalized periodic boundary value problems}

 Consider the following problem (3.1)(3.6)
     \begin{eqnarray}
  &&(\Lambda(t)x')'+B(t)x=0\nn\\
  &&x(1)=Mx(0), x'(1)=Nx'(0) \lb{3.6}
    \end{eqnarray}
where $M\in GL(n), M^{\tau}\Lambda(1)N=\Lambda(0),\Lambda\in
C([0,1];GL_s(n))$ and $\Lambda(t)$ is positive definite. Define
$X:=L^2((0,1);\R^n), Y:=\{x:[0,1]\to\R^n|(\Lambda(t)x'(t))'\in
L^2(0,1;\R^n)$ and $x$ satisfies $(3.6)\}$. Then $Y\hookrightarrow
X$ is compact. Define $(Ax)(t):=(\Lambda(t)x'(t))'$ for every
$x\in Y$. Then $A:Y\to X$ is continuous and
$X=\ker(A)\bigoplus\Im(A)$. A simple calculation shows that
 \be
 \int_0^1[(\Lambda(t)x'(t))'\cdot x(t)]dt\leq 0\lb{3.7}
 \ee
for every $x\in Y$. So similar to definition 3.2 we have from proposition 1.4 the following
definition.

{\bf Definition 3.6} For any $B \in L^{\infty}(0,1;GL_s(\R^n))$ we
define
   \begin{eqnarray}
  &&\nu^s_{\Lambda,M}(B):=\dim Ker(A+B),\nn\\
  &&i^s_{\Lambda,M}(B):=\sum_{\lambda<0}\nu^s_{\Lambda,M}(B+I_n\lambda).\nn
  \end{eqnarray}

{\bf Proposition 3.7}. (1)For any $B\in
L^{\infty}((0,1);$GL$_s(\R^n))$, we have that $E^0(\Lambda,B,M)$
is the solution subspace of (3.1)(3.6) and
$\nu^s_{\Lambda,M}(B)\in\{0,1,2,\cdots,2n\}$.

(2) For any $B_1, B_2\in L^{\infty}((0,1);$GL$_s(\R^n))$, if
$B_1\leq B_2$, we have $i^s_{\Lambda,M}(B_1)\leq
i^s_{\Lambda,M}(B_2)$; if $B_1<B_2$, we have
$i^s_{\Lambda,M}(B_1)+\nu^s_{\Lambda,M}(B_1)\leq
i^s_{\Lambda,M}(B_2)$.

 (3)For any $\Lambda_1,\Lambda_2$ with
$\Lambda_1(1)=\Lambda_1(0),\Lambda_2(1)=\Lambda_2(0)$, if
$\Lambda_1<\Lambda_2$, then
$i^s_{\Lambda_1,M}(B)+\nu^s_{\Lambda_1,M}(B)\leq
i^s_{\Lambda_2,M}(B)$.

(4)If $B_i\in L^{\infty}((0,1);$GL$_s(\R^{n_i})), \Lambda_i\in
C([0,1];$GL$_s(\R^{n_i})),M_i,N_i\in $GL$(\R^{n_i}))$ with
$M_i^T\Lambda_i(1)N_i=\Lambda_i(0), i=1,2$ and
$B=\diag\{B_1,B_2\},\Lambda=\diag\{\Lambda_1,\Lambda_2\},M=\diag\{M_1,M_2\},N=\diag\{N_1,N_2\}$
then
$i^s_{\Lambda,M}(B)=i^s_{\Lambda_1,M_1}(B_1)+i^s_{\Lambda_2,M_2}(B_2),
\nu^s_{\Lambda,M}(B)=\nu^s_{\Lambda_1,M_1}(B_1)+\nu^s_{\Lambda_2,M_2}(B_2)$.

{\bf Example 3.8}. Let
$\alpha_1\leq\alpha_2\leq\cdots\leq\alpha_n$ be the eigenvalues of
a constant matrix $A$. Then
   \begin{eqnarray}
  &&i^s_{\lambda I_n,I_n}(A)={^\#}\{k:\alpha_k>0\}+2\sum_{k=1}^n{^\#}\{j\in\N:4\lambda j^2\pi^2<\alpha_k\},\nn\\
  &&\nu^s_{\lambda I_n,I_n}(A)={^\#}\{k:\alpha_k=0\}+2\sum_{k=1}^n{^\#}\{j\in\N:4\lambda j^2\pi^2=\alpha_k\},\nn\\
  &&i^s_{\lambda I_n,-I_n}(A)=2\sum_{k=1}^n{^\#}\{j\in\N:\lambda(2j-1)^2\pi^2<\alpha_k\},\nn\\
  &&\nu^s_{\lambda I_n,-I_n}(A)=2\sum_{k=1}^n{^\#}\{j\in\N:\lambda(2j-1)^2\pi^2=\alpha_k\}.\nn
    \end{eqnarray}
where ${^\#}S$ denotes the number of elements in a set $S$. For
$a\in\R\setminus\{\pm1,0\}$, we have with
$\mu_0=\arccos{2\over{a^{-1}+a}}$ that
  \begin{eqnarray}
  &&i^s_{\lambda I_n,a I_n}(A)=\sum_1^k{^\#}\{j\in
  \N:\lambda(2j\pi+\mu_0)^2<\alpha_k\}+\sum_{k=1}^n{^\#}\{j\in\N:\lambda(2\pi-\mu_0+2j\pi)^2<\alpha_k\},\nn\\
  &&\mu^s_{\lambda I_n,a I_n}(A)=\sum_1^k{^\#}\{j\in\N:\lambda(2j\pi+\mu_0)^2=\alpha_k\}
  +\sum_{k=1}^n{^\#}\{j\in\N:\lambda(2\pi-\mu_0+2j\pi)^2=\alpha_k\}.\nn
  \end{eqnarray}
{\bf Remark 3.9}. The first two formulae in Example 3.8 were given
first by Mawhin and Willem in the book[43] when $\lambda=1$. In
order to discuss minimal periodic solution problems Y. Long[44,45]
established two kind of index theory for linear Hamiltonian
systems satisfying periodic boundary value conditions in some
sense of symmetries in 1993 and 1994. \vskip 1.5cm

We discuss solvability of the following nonlinear systems
(3.8)(3.6):
 \begin{eqnarray}
  &&(\Lambda(t)x')'+B(t,x)x+h(t,x)=0,\lb{3.8}\\
  &&x(1)=Mx(0),x'(1)=Nx'(0)\nn
  \end{eqnarray}
where $A:\in C([0,1]\times\R^n,GL_s(\R^n)), h:[0,1]\times \R^n\to
\R^n$ are continuous. Generally, (3.7)(3.6) is not a Lagrangian
system, i.e., we can not find a  $V\in C^1([0,1]\times\R^n,\R)$
such that $V'(t,x)=+A(t,x)x+h(t,x) $. Even though we still have
the following theorem, which proof is similar to theorem 1.6's.

{\bf Theorem 3.10} Assume

  (1) there exist $B_1, B_2\in
L^{\infty}((0,1);$GL$_s(\R^n))$ with $B_1\leq B_2,
i^s_{\Lambda,M}(B_1)=i^s_{\Lambda,M}(B_2)$,
$\nu^s_{\Lambda,M}(B_2)=0$ such that
  $$
  B_1(t)\leq B(t,x)\leq B_2(t), x\in R^n,\,\,\hbox{a.e.}\,\, t\in (0,1);
  $$
  (2) $h(t,x)=(|x|)$ as $|x|\to+\infty$.
Then (3.8)(3.6) has at least one solution.

{\bf Example 3.11} Let
$B(t,x)=B_1(t)\cos^2|x|^2+B_2(t)\sin^2|x|^2,
h(t,x)=x(1+|x|^2)\sin|x|t$. As $\Lambda(t)=I_n,M=N=-I_n$, choose
$B_1(t)=(\pi^2(2k-1)^2+\epsilon)I_n,
B_2(t)=(\pi^2(2k+1)^2-\epsilon)I_n$; as $\Lambda(t)=I_n,M=N=I_n$,
choose $B_1(t)=(4\pi^2k^2+\epsilon)I_n,
B_2(t)=(4\pi^2(k+1)^2-\epsilon)I_n$; as $\Lambda(t)=I_n,M=\lambda
I_n, N=(\lambda)^{-1}I_n$ with $\lambda\in\R\setminus\{\pm1,0\}$,
choose $B_1(t)=((2k\pi+\mu)^2+\epsilon)I_n,
B_2(t)=((2k\pi+2\pi-\mu)^2-\epsilon)I_n$ with
$\mu=\arccos{2\over{\lambda^{-1}+\lambda}}$. Then (3.8)(3.6) has
at least one solution provided $\epsilon>0$ is sufficiently small.

Finally, we will consider the following Lagrangian system
(3.5)(3.6)
 \begin{eqnarray}
  &&(\Lambda(t)x')'+V'(t,x)=0,\nn\\
  &&x(1)=M_1x(0),x'(1)=M_2x'(0)\nn
  \end{eqnarray}

From theorem 1.7 and its proof we have the following theorem.

{\bf Theorem 3.12} Assume

(1) $V\in C^2([0,1]\times \R^n, \R), B_1(t)\leq V''(t,x)\leq
B_2(t) $ for $|x|\geq r>0$ with
$i^s_{\Lambda,M}(B_1)=i^s_{\Lambda,M}(B_2),\nu^s_{\Lambda,M}(B_2)=0$.

(2) $V'(t,0)\equiv 0, {\bar B}(t):=V''(t,0)$ and
$i^s_{\Lambda,M}(B_1)\notin [i^s_{\Lambda,M}({\bar B}),
i^s_{\Lambda,M_1}({\bar B})+\nu^s_{\Lambda,M_1}({\bar B})]$.

Then problem (3.5)(3.6) has at least one nontrivial solution.
Moreover, if we assume

(3)$\nu^s_{\Lambda,M}({\bar B})=0,|i^s_{\Lambda,M}({\bar
B})-i^s_{\Lambda,M}({\bar B})|\geq 2n$.

Then (3.5)(3.6) has two nontrivial solutions.

{\bf Remark 3.13}. I. Ekeland[5 suggested to discuss the following
boundary value conditions
    \be
\left(\matrix{x(1)\cr x'(1)\cr}\right)=M\left(\matrix{x(0)\cr
x'(0)\cr}\right),\lb{3.9}
    \ee
where  $M\in GL(2n)$ satisfying
 $$
M^T\left(\matrix{0&-\Lambda(1)\cr
\Lambda(1)&0\cr}\right)M=\left(\matrix{0&-\Lambda(0)\cr
\Lambda(0)&0\cr}\right)
 $$
Condition (3.6) is a special case of (3.9). This condition is
chosen because we can get an inequality like (3.7), and so we can
establish an index theory like in definition 3.6. In next section
we will discuss generalized periodic boundary condition for first
order Hamiltonian system, which will cover condition (3.9).

\setcounter{equation}{0}
\section{First order Hamiltonian systems}

{\bf 4.1 Bolza  BVPs}

 In this subsection we will establish
  a classification theory for the following Hamiltonian system
    \begin{eqnarray}
  &&\dot{x}=JB(t)x\lb{4.1}\\
  &&x_1(0)\cos\alpha+x_2(0)\sin\alpha=0\lb{4.2}\\
  &&x_1(1)\cos\beta+x_2(1)\sin\beta=0\lb{4.3}
 \end{eqnarray}
where $B\in L^{\infty}((0,1);GL_s(\R^{2n})), 0\leq\alpha<\pi$ and
$0<\beta\leq\pi, x=(x_1,x_2)\in\R^n\times\R^n$. Let
$X:=L^2((0,1);\R^{2n}), Y=\{x:[0,1]\to\R^{2n}|x'\in
L^2((0,1);\R^{2n}),x(t)$ satisfies $(4.2)(4.3)\}$. Define $A:Y\to
X$ by $(Ax)(t):=Jx'(t)$. We can choose suitable value $\lambda\in
\R$ such that problem (4.1)(4.2)(4.3) with $B(t)$ replaced by
$\lambda I_{2n}$ has no nontrivial solutions. From general theory
of ordinary differential equations, for any $h\in
L^2((0,1),\R^{2n})$ the following problem (4.2)(4.3) and
  $$
  J\dot{x}+\lambda x=h(t)
  $$
has a unique solution. So from lemma 2.1 $A$ is continuous and closed and,
$\ker(A)\bigoplus\Im(A)=X$.

 {\bf Definition 4.1} For any $B\in L^{\infty}((0,1);GL_s(\R^{2n}))$,
 we define
  \begin{eqnarray}
  &&\nu^f_{\alpha,\beta}(B):=\dim \ker(A+B),\nn\\
  &&i^f_{\alpha,\beta}(\diag\{0,I_n\}):=i^s_{I_n,\alpha,\beta}(0),\nn\\
&&i^f_{\alpha,\beta}(B):=i^f_{\alpha,\beta}(\diag\{0,I_n\})+I^f_{\alpha,\beta}(\diag\{0,I_n\},B);\nn
  \end{eqnarray}
and
  \begin{eqnarray}
  &&I^f_{\alpha,\beta}(B_1,B_2)=\sum_{\lambda\in [0,1)}\nu^f_{\alpha,\beta}((1-\lambda)B_1+\lambda B_2)
  \,\,\,\,\hbox{as}\,\,\,\,B_1<B_2,\nn\\
  &&I^f_{\alpha,\beta}(B_1,B_2)=I^f_{\alpha,\beta}(B_1,k id)-I^f_{\alpha,\beta}(B_2,k id)
          \,\,\,\,\hbox{for every}\,\,\,\, B_1,B_2\,\,\,\,\hbox{with}\,\,\,\,k id >B_1, k id>B_2.\nn
  \end{eqnarray}

{\bf Proposition 4.2} We have the following property:

 (1) For any $B\in L^{\infty}((0,1);GL_s(\R^{2n}))$, we have
  $$
 \nu^f_{\alpha,\beta}(B)\in\{0,1,\cdots,n\}.
  $$

 (2) For any $B_1, B_2\in L^{\infty}((0,1);$GL$_s(\R^{2n}))$
satisfying $B_1< B_2$, we have
$i^f_{\alpha,\beta}(B_1)+\nu^f_{\alpha,\beta}(B_1)\leq
i^f_{\alpha,\beta}(B_2)$.

 (3) For any $B\in L^{\infty}((0,1);$GL$_s(\R^n))$, we have
   $$
   (i^f_{\alpha,\beta}(\diag\{B,I_n\}),\nu^f_{\alpha,\beta}(\diag\{B,I_n\}))
   =(i^s_{I_n,\alpha,\beta}(B),\nu^s_{I_n,\alpha,\beta}(B)).
   $$

{\bf Proof}. We only prove
  \be
  i^f_{\alpha,\beta}(\diag\{B,I_n\})
   =i^s_{I_n,\alpha,\beta}(B).\lb{4.4}
   \ee

Case 1: $B>0$. Choose a negative number $c\in\R$. Similar to the
(iii) of theorem 2.4, we have
  \begin{eqnarray}
&&i^f_{\alpha,\beta}(\diag\{B,I_n\})-i^f_{\alpha,\beta}(\diag\{0,I_n\})\nn\\
&&=i^f_{\alpha,\beta}(\diag\{B,I_n\}|cI_{2n})-i^f_{\alpha,\beta}(\diag\{0,I_n\}|cI_{2n})\nn\\
&&=\sum_{\lambda\in[0,1)}\nu^f_{\alpha,\beta}(\diag\{\lambda
       B,I_n\})\nn\\
&& =\sum_{\lambda\in[0,1)}\nu^s_{I_n,\alpha,\beta}(\lambda
        B)\nn\\
&&=i^s_{I_n,\alpha,\beta}(B)-i^s_{I_n,\alpha,\beta}(0)\nn
  \end{eqnarray}
Combining the second formula in definition 4.1, formula(4.4)
follows.

Case 2: $B$ is arbitrary. Choose a positive number $c$ such that
$cI_n>B$. Similar to Case 1, we have
   \begin{eqnarray}
&&i^f_{\alpha,\beta}(\diag\{c_1 I_n,I_n\})-i^f_{\alpha,\beta}(\diag\{B,I_n\})\nn\\
&&=i^s_{I_n,\alpha,\beta}(c_1I_n)-i^s_{I_n,\alpha,\beta}(B)\nn
  \end{eqnarray}
and hence formula (4.4). \hfill\hb

  We now begin to discuss solvability of the following nonlinear
Hamiltonian systems:
 \begin{eqnarray}
  &&x''=JH'(t,x),\lb{4.4}\\
  &&x_1(0)\cos\alpha+x_2(0)\sin\alpha=0\nn\\
  &&x_1(1)\cos\beta+x_2(1)\sin\beta=0\nn
  \end{eqnarray}
where $H:[0,1]\times\R^{2n}\to \R^{2n}$ is differentiable and $H'(t,x)$ is the gradient of $H$
with respect to $x$.

From theorems 1.6 1nd 1.7 we have the following theorems.

 {\bf Theorem 4.3} Assume

  (1) there exist $B_1, B_2\in
L^{\infty}((0,1);$GL$_s(\R^{2n}))$ with $B_1\leq B_2,
i^f_{\alpha,\beta}(B_1)=i^f_{\alpha,\beta}(B_2)$,
$\nu^f_{\alpha,\beta}(B_2)=0$ such that

  (2) $H'(t,x)-B(t,x)$ is bounded, $B:[0,1]\times\R^{2n}\to
  GL_s(\R^{2n})$ is continuous and
  $$
  B_1(t)\leq B(t,x)\leq B_2(t), (t,x)\in (0,1)\times\R^{2n},\,\hbox{with}\,\, |x|\geq r>0.
  $$
  Then (4.4)(4.2)(4.3) has at least one solution.

{\bf Theorem 4.4}  Assume

(1) $H\in C^2([0,1]\times \R^{2n}, \R), B_1(t)\leq H''(t,x)\leq
B_2(t) $ for $|x|\geq r>0$ with
$i^f_{\alpha,\beta}(B_1)=i^f_{\alpha,\beta}(B_2),\nu^f_{\alpha,\beta}(B_2)=0$.

(2) $H'(t,0)\equiv 0, {\bar B}(t):=H''(t,0)$ and
$i^f_{\alpha,\beta}(B_1)\notin [i^f_{\alpha,\beta}({\bar B}),
i^f_{\alpha,\beta}({\bar B})+\nu^f_{\alpha,\beta}({\bar B})]$.

Then problem (4.4)(4.2)(4.3) has at least one nontrivial solution.
Moreover, if we further assume

(3)$\nu^f_{\alpha,\beta}({\bar
    B})=0,|i^f_{\alpha,\beta}(B_1)-i^f_{\alpha,\beta}({\bar B})|\geq
     n$.

Then (4.4)(4.2)(4.3) has two nontrivial solutions.

Note that in [7] we discussed the special case $\alpha=0,\beta=\pi$.

{\bf 4.2 Generalized periodic boundary value problems}

   Consider the following problem (4.1)(4.6)
     \begin{eqnarray}
  &&x'=JB(t)x\nn\\
  &&x(1)=Px(0) \lb{4.5}
    \end{eqnarray}
where $P\in Sp(2n)$ is prescribed. Define $X:=L^2((0,1);\R^{2n}),
Y:=\{x:[0,1]\to\R^{2n}|x'\in L^2(0,1;\R^{2n})$ and $x$ satisfies
$(4.6)\}$. Then $Y\hookrightarrow X$ is compact. Define
$(Ax)(t):=Jx'(t)$ for every $x\in Y$. Similar to Proposition 7 in
page 22 of Ekeland's book[5], for the given $P\in Sp(2n)$ there
exists $\lambda\in\R$ such that $e^{J\lambda}-P$ is invertible. So
(4.1)(4.6) with $B(t)$ replaced by $\lambda I_{2n}$ has only the
trivial solution. Thus, from lemma 2.1 $A:Y\to X$ is continuous,
closed and $X=\ker(A)\bigoplus\Im(A)$.

Choose $i^f_P(0):=i_P(I_{2n})$ defined by definition 2.2 in [8].
we have the following definition.

{\bf Definition 4.5} For any $B\in
L^{\infty}((0,1);$GL$_s(\R^{2n}))$, we define
  \begin{eqnarray}
  &&\nu^f_P(B)=dim \ker(A+B),\nn\\
  &&i^f_P(B)=i^f_P(0)+I^f_P(0,B);\nn
  \end{eqnarray}
and
  \begin{eqnarray}
  &&I^f_P(B_1,B_2)=\sum_{\lambda\in [0,1)}\nu^f_P((1-\lambda)B_1+\lambda B_2)
  \,\,\,\,\hbox{as}\,\,\,\,B_1<B_2,\nn\\
  &&I^f_P(B_1,B_2)=I^f_P(B_1,k id)-I^f_P(B_2,k id)
          \,\,\,\,\hbox{for every}\,\,\,\,B_1,B_2\,\,\,\,\hbox{with}\,\,\,\, k id >B_1, k id>B_2.\nn
  \end{eqnarray}
From theorem 1.5 we have the following proposition.

{\bf Proposition 4.6}. (1)For any $B\in
L^{\infty}((0,1);$GL$_s(\R^{2n}))$, we have
$\nu^f_P(B)\in\{0,1,2,\cdots,2n\}$.

(2) For any $B_1, B_2\in L^{\infty}((0,1);$GL$_s(\R^{2n}))$
satisfying $B_1< B_2$, we have $i^f_P(B_1)+\nu^f_P(B_1)\leq
i^f_P(B_2)$.

We now discuss solvability of the following nonlinear system
(4.6)(4.5):
 \begin{eqnarray}
  &&x'=JH'(t,x)\nn\\
  &&x(1)=Px(0)\nn
  \end{eqnarray}
where $H:[0,1]\times \R^{2n}\to \R$ is differentiable and $P\in
Sp(2n)$ is prescribed.

Similar to theorems 4.3 and 4.4 we have

{\bf Theorem 4.7} Assume

  (1) there exist $B_1, B_2\in
L^{\infty}((0,1);$GL$_s(\R^{2n}))$ with $B_1\leq B_2,
i^f_P(B_1)=i^f_P(B_2)$, $\nu^f_P(B_2)=0$ such that
  $$
  B_1(t)\leq B(t,x)\leq B_2(t), x\in R^n,\,\,\hbox{a.e.}\,\, t\in (0,1);
  $$
  (2) $H'(t,x)-B(t,x)$ is bounded.

Then (4.4)(4.5) has at least one solution.

{\bf Theorem 4.8}  Assume

(1) $H\in C^2([0,1]\times \R^{2n}, \R), B_1(t)\leq H''(t,x)\leq
B_2(t) $ for $|x|\geq M>0$ with
$i^f_P(B_1)=i^f_P(B_2),\nu^f_P(B_2)=0$.

(2) $H'(t,0)\equiv 0, {\bar B}(t):=H''(t,0)$ and $i^f_P(B_1)\notin
[i^f_P({\bar B}), i^f_P({\bar B})+\nu^f_P({\bar B})]$.

Then problem (4.4)(4.5) has at least one nontrivial solution.
Moreover, if we assume

(3)$\nu^f_P({\bar B})=0,|i^f_P(B_1)-i^f_P({\bar B})|\geq 2n$.

Then (4.4)(4.5) has two nontrivial solutions.

Note that in [8] the above problem have been discussed already separately.

\setcounter{equation}{0}
\section{Second order elliptic partial differential equations}

In this section we will discuss index theory for linear elliptic
equations satisfying Dirichlet boundary conditions and nontrivial
solutions for nonlinear elliptic equations. First we consider the
following linear systems:
    \begin{eqnarray}
  &&\Delta u+b(x)u=0,x\in \Omega\lb{5.1}\\
  &&u|{\partial \Omega}=0\lb{5.2}
  \end{eqnarray}
where $\Omega\in \R^n$ is a bounded open domain, and its boundary
$\partial\Omega$ is smooth, $b\in L^{\infty}(\Omega)$.

Define $X:=L^2(\Omega), Y:=H^2_0(\Omega)$ and $Au=\triangle u,
(Bu)(x)=b(x)u(x)$. Then the embedding $Y\hookrightarrow X$ is
compact, $A:Y\to X$ and $B:X\to X$ are continuous and
self-adjoint. It is well-known that the spectrum
$\sigma(-A)\subset(0,+\infty)$. And by the Dirichlet Principle,
for any $f\in L^2(\Omega)$ equation $\Delta u=f$ and (5.2) has a
weak solution. This weak solution is also classical solution. So
$Im(A)=X$ and $\ker(A)=\{\theta\}$. From proposition 1.4 we can
give the following definition.

{\bf Definition 5.1} For any $b\in L^{\infty}(\Omega)$, we define
    \begin{eqnarray}
  &&\nu_{\Delta}(b)=\dim\ker(A+B)\nn\\
  &&i_{\Delta}(b)=\sum_{\lambda<0}\nu_{\Delta}(b+\lambda).\nn
  \end{eqnarray}

The following proposition comes from proposition 1.5 directly.
Note that for any $b_1,b_2\in L^{\infty}(\Omega)$, we define
$b_1\leq b_2$ if and only if $b_1(x)\leq b_2(x)$ for a.e. $x\in
\Omega$; and define $b_1<b_2$ if and only if $b_1\leq b_2$ and
$b_1(x)<b_2(x)$ on a subset of $\Omega$ with positive measure.

{\bf Proposition 5.2} (1) $\nu_{\Delta}(b)$ is finite.

(2)For any $b_1<b_2$ belonging to $L^{\infty}(\Omega)$, we have
  $$
  \nu_{\Delta}(b_1)+i_{\Delta}(b_1)\leq i_{\Delta}(b_2).
  $$

Finally, we consider the following problem:
  \begin{eqnarray}
  &&\Delta u+f(x,u)=0,x\in \Omega\lb{5.3}\\
  &&u|{\partial \Omega}=0\nn
  \end{eqnarray}

Define $F(x,u):=\int_0^uf(x,s)ds$ and
$\Phi(u)=\int_{\Omega}F(x,u(x))dx$. Then $\Phi'(u)=f(\cdot,u(\cdot)),
\Phi''(u)={\partial\over\partial u}f(\cdot,u(\cdot))$ and equation (5.3)(5.2) is
equivalent to (1.2).

 {\bf Theorem 5.3} Assume that
  $$
  b_1(x)\leq f(x,u)/u\leq b_2(x), |u|>r>0
  $$
and $i_{\Delta}(b_1)=i_{\Delta}(b_2), \nu_{\Delta}(b_2)=0$. Then
(5.3)(5.2) has at least one solution.

 {\bf Proof}. Define $g(x,u)=f(x,u)/u$ as $|u|>r$; $g(x,u)=b_1(x)x$ as $|u|\leq
 r$, and
 $B(u)=g(\cdot,u(\cdot))$. Then from theorem 1.6 and its proof we can complete the proof.\hfill\hb

And from theorems 1.7, 1.8 and 1.9 we obtain the following
results.

   {\bf Theorem 5.4}  Assume

(1) $f\in C^1(\Omega\times \R, \R), b_1(x)\leq {\partial\over\partial u}f(x,u)\leq b_2(x)
$ for $|u|\geq r>0$ with
$i_{\Delta}(b_1)=i_{\Delta}(b_2),\nu_{\Delta}(b_2)=0$.

(2) $f(x,0)\equiv 0, {\bar b}(x):={\partial\over\partial u}f(x,0)$ and
$i_{\Delta}(b_1)\notin [i_{\Delta}({\bar b}), i_{\Delta}({\bar
b})+\nu_{\Delta}({\bar b})]$.

Then (5.3)(5.2) has at least one nontrivial solution.

{\bf Theorem 5.5} Assume that

(1) $f'\in C(\Omega\times \R,\R)$ and there exist $b_1, b_2\in
L^{\infty}(\Omega)$ with $\nu_{\Delta}(b_1)=0$ such that
   $$
   b_1\leq f;(x,u)\leq b_2   \forall (x,u)\in \Omega\times \R;
   $$

 (2)there exists $b_3\in L^{\infty}(\Omega)$ with $b_1<b_3$
 and $i_{\Delta}(b_1)=i_{}(b_3), \nu_{\Delta}(b_3)=0$ such that
  $$
  \Phi(x)\leq {1\over 2}(b_3(x)x,x)+c \forall x\in X;
   $$
  (3) $f(x,0)=0,
  {\partial\over\partial u}f(x,0)>b_1(x),\nu_A({\partial\over\partial u}f(\cdot,0))=0$
  and
  $i_{\Delta}({\partial\over\partial u}f(\cdot,0))>i_{\Delta}(b_1)$.

  Then (5.3)(5.2) has two distinct nontrivial solutions.

{\bf Theorem 5.6} Assume that

  (1) there exist  $b_1,b_2\in L^{\infty}(\Omega)$ satisfying
 $b_1\leq b_2$ and
 $i_{\Delta}(b_1)+\nu_{\Delta}(b_1)=i_{\Delta}(b_2), \nu_{\Delta}(b_2)=0$ such that $\int_0^uf(x,s)ds-{1\over
 2}(b_1(x)u^2$ is convex with respect to u and
   $$
    \int_0^uf(x,s)ds\leq {1\over 2}b_2(x)u^2+c \forall (x,u)\in \Omega\times\R.
   $$
Then (5.3)(5.2) has a solution.

   Moreover, if we further assume that

   (2) $f(x,0)=0$ and there exists $b_0\in L^{\infty}(\Omega)$ satisfying
$b_0\geq b_1$ and
   $$
   i_{\Delta}(b_0)>i_{\Delta}(b_1)+\nu_{\Delta}(b_1).
   $$
Then (5.3)(5.2) has at least one nontrivial solution.

{\bf Remark 5.7} Theorems 5.4, 5.5 and 5.6 cover some results in
[3, Chapter III].

{\bf Acknowledgement} Part of the manuscript was finished during
my stay at IHES and Universite Paris-Dauphine from Aug 2004 to Oct
2005. I would like to express my sincere thanks to Profs
Jean-Pierre BOURGUIGNON and Eric Sere and other members and
visitors of the two institutes for their warm helps. Special
thanks are devoted to Eric Sere for offering a proof for lemma
3.1. During the preparation of the paper I also visited Chern
Institute of Mathematics invited by Prof Yiming Long. I also would
like to express my thanks to Prof Yiming Long and Weiping Zhang
for their hospitality. The final manuscript was finished while my visiting
at PIMS. I thank Prof Ivar Ekeland for his invitation and help.

{\bibliographystyle{abbrv}

\begin{thebibliography}{99}

\bibitem{AmanZ} H. Amann and E. Zehnder, Nontrivial soultions for
a class of nonresonance problems and applications to nonlinear
differential equations, Annali Scuola Norm. Sup. Pisa
7(1980)439-603.

\bibitem{Chang1} K. C. Chang, Solutions of asymptotically linear
operator equations via Morse theory, Comm. Pure Appl. Math.
34(1981)693-712.

\bibitem{Chang2} K. C. Chang, Infinite dimentional Morse theory
and multiple solution problems. Birkhauser. Basel(1993).

\bibitem{Lo5} K. C. Chang, Critical point theory and its application,
Shanghai Sci. Tech. Press(1986)(in Chinese).


\bibitem{Eke1}I. Ekeland, Convexity methods in Hamiltonian
mechanics. Springer-Verlag. Berlin. 1990.


\bibitem{Do1} Y. Dong, Index theory, nontrivial solutions and
asymptotically linear second order Hamiltonian systems. J. Differ.
Equations 214(2005)233-255.


\bibitem{Do2} Y. Dong, Maslov type index theory for linear
Hamiltonian systems with Bolza boundary value conditions and
multiple solutions for nonlinear Hamiltonian systems. Pacific J
Math (2005)253-280.

\bibitem{Do3} Y. Dong, $P-$index theory for linear Hamiltonian
systems and multiple solutions for nonlinear Hamiltonian systems.
Nonlinearity 19(2006)1275-1294.


\bibitem{Eke3}I. Ekeland, Une theorie de Morse pour les systemes
hamiltoniens convexes. Ann IHP "Analyse non lineaire"
1(1984)19-78.

\bibitem{ColZehn} C. Conley and E Zehnder, Morse-type index theory
for flows and periodic solutions for Hamiltonian equations. Comm.
Pure Appl. Math. 37(1984)207-253.

\bibitem{LoZ} Y. Long and E. Zehnder, Morse theory for forced oscillations of
asymptotically linear Hamiltonian systems. Stock. Process. Phys.
Geom. ed S Alberverio et al(Teaneck, NJ:World
Scientific)(1990)528-563.

\bibitem{Long} Y. Long, Maslov-type index, degenerate critical
poin ts, and asymptotically linear Hamiltonian systems. Sci. China
33(1990)1409-1419.

\bibitem{Long} Y. Long, A Maslov-type index theory for symplectic
p aths. Topol. Methods Nonlinear Anal. 10(1997)47-78.

\bibitem{Ch2}I. Ekeland and H. Hofer, Periodic solutions with prescribed period
for convex autonomous Hamitonian systems Invent. Math.
81(1985)155-188.

\bibitem{Ch2} D. Dong and Y. Long, The iteration formula of
Maslov-type index theory with applications to nonlinear
Hamiltonian systems Trans. American Math. Soc. 349(1997)2619-2661.

\bibitem{Ch2} I. Ekeland and H. Hofer, Convex Hamiltonian energy
surfaces and their closed trajectories. Comm. Math. Phys.
113(1987)419-467.

\bibitem{Ch2} Y. Long and C. Zhu, Closed characteristics on
compact  convex hypersurfaces in $\R^{2n}$. Ann. Math.
155(2002)317-368.

\bibitem{Ch2} C. Liu, Y. Long and C. Zhu, Multiplicity of closed
characteristics on symmetric convex hypersurfaces in $\R^{2n}$.
Math. Ann. 323(2002)201-215.


\bibitem{FEI} G. Fei, Relative Morse index and its applications to the
Hamiltonian systems in the presnece of symmetry. J. Diff. Equa.
122(1995)302-315.


\bibitem{FEI} G. Fei, Maslov-type index and periodic solution of
asymptotically linear Hamiltonian systems which are resonant at
infinity, J. Differential Equations 121(1995)121-133.


\bibitem{Ch2} J. Su, Nontrivial periodic solutions for the
asymptotically linear Hamiltonian systems with resonance at
infinity. J. Differential Equations 145(1998)252-273.

\bibitem{Guo} Y. Guo, Nontrivial periodic solutions for
asymptotically linear Hamiltonian systems with resonance. J.
Differential Equations 175(2001)71-87.

\bibitem{Lo} Y. Long, Index theory for symplectic paths with
applications, Progress in Math. No. 207, Birkh\"auser. Basel.
2002.


\bibitem{ZL} C. Zhu and Y. Long, Maslov type index theorey for symplectiuc paths and
spectral flow(I). Chinese Ann. of Math. 20B(1999)413-424.

\bibitem{LZ} Y. Long and C.Zhu, Maslov type index theorey for symplectiuc paths and
spectral flow(II). Chinese Ann. of Math. 21B(2000) 89-108.

\bibitem{Ch2} S. Cappell, Lee R and Miller E Y, On the Maslov
index. Comm. Pure Appl. Math. 17(1994)121-186.

\bibitem{Ch2} J. Leray, Lagrangian Analysis and quantum mechanics, a
mathematical structure related to asymptotic expansions and the
Maslov index(Cambridge, MA:MIT Press)1981.

\bibitem{Ch2} P. Dazord, Invariants homotopiques attachs aux fibrs
symplectiques Ann. Inst. Fourier 29(1979)25-78.

\bibitem{Ch2} de Gosson M, The structure of $q-$symplectic
geometry J. Math Pures Appl. 71(1992)429-453.

\bibitem{Hartman} P. Hartman, Ordinary differential equations.
Second edition(1982). Birkhauser. Boston Basel Stuttgart.


\bibitem{Eke2} I. Ekeland,  N. Ghoussoub and H. Tehrani,
Multiple solutions for a classical problem in the calculus of
variations, J. Differential Equations 131 (1996)229-243

\bibitem{Eke3} F. Clarke and E. Ekeland, Nonlinear oscillations
and boundary value problems for Hamiltonian systems, Arch.
Rational Mech. Anal. 78(1982)315-333.


\bibitem{EL} Z. Wang, Multiple solutions for infinite functional
and applications to asymptotically linear problems, Math.
Sinica(N.S.)5(1989)101-113.


\bibitem{Dong}Y. Dong, On Equivalent Conditions for the Solvability
 of Equation $(p(t)x')'+f(t,x)=h(t)$ Satisfying
Linear Boundary Conditions with f Restricted by Linear Growth
Conditions, J. Math. Anal. Appl. 245 (2000)204-220.

\bibitem{Dong} Y. Dong,  On the solvability of asymptotically
positively homogeneous equations with S-L boundary value
conditions, Nonlinear Analysis 42(2000) 1351-1363.


\bibitem{Hartman}H. Wang and Y. Li, Existence and uniqueness of periodic
     solutions for Duffing equations across many points of resonance.
     J. Differential Equations 108 (1994)152--169


\bibitem{wangl}H. Wang and Y. Li, Two-point boundary value problems for
     second order ordinary differential equations across many resonant points.
     J. Math. Anal. Appl. 179 (1993)61--75


\bibitem{Fa}C. Fabry, Landesman-Lazer conditions for periodic
     boundary value problems with asymmetric nonlinearities.
     J. Differential Equations 116 (1995)405--418


\bibitem{Ville}S. Villegas, A Neumann problem with asymmetric
     nonlinearity and a related minimizing problem.
     J. Differential Equations 145 (1998)145--155


\bibitem{IaNka}R. Iannacci, M. Nkashama, Nonlinear elliptic
     partial differential equations at resonance:
     higher eigenvalues. Nonlinear Anal. 25 (1995)455--471

\bibitem{IaNka}R. Iannacci, M. Nkashama, and  J. Ward, Nonlinear
      second order elliptic partial differential equations at resonance.
      Trans. Amer. Math. Soc. 311 (1989)711--726


\bibitem{NkaR}M. Nkashama, S. Robinson, Resonance and nonresonance in terms of average values.
     J. Differential Equations 132 (1996)46--65


\bibitem{Lo6} J. Mawhin and M. Willem, Critical point theory and
Hamiltonian systems. Springer. Berlin 1998.


\bibitem{Lo6} Y. Long, The minimal period problem for classical
Hamiltonian systems with even potentials. Ann. Inst. H. Poincare
Anal. non lineaire. 10(1993)605-626.

\bibitem{Lo6} Y. Long, The minimal period problem of periodic
solutions for autonomous superquadratic second order Hamiltonian
systems. J. Differential Equations 111(1994)147-174.









\end{thebibliography}

\end{document}